\documentclass[12pt,a4paper,reqno]{amsart}
\usepackage{amsmath,amsthm}
\usepackage{amsfonts}
\usepackage{amssymb}
\usepackage{tikz-cd}
\allowdisplaybreaks

\newcommand{\be}{\begin{equation}}
\newcommand{\ef}{\end{equation}}
\chardef\bslash=`\\ 

\newtheorem*{thm*}{Theorem}

\theoremstyle{definition}

\newtheorem*{remark*}{Remarks}
\newtheorem*{defn*}{Definition}

\theoremstyle{remark}

\newcommand{\wt}{\widetilde}
\newcommand{\wh}{\widehat}
\newcommand{\fc}{\frac}

\newcommand{\iy}{\infty}
 \renewcommand{\sectionmark}[1]{}
\renewcommand{\Im}{\operatorname{Im}}

\newcommand{\const}{\operatorname{const}}

 \usepackage{amsfonts}

\newcommand{\D}{\mathbb{D}}
\newcommand{\om}{\omega}
\newcommand{\z}{\zeta}
\newcommand{\ov}{\overline}
\newcommand{\vp}{\varphi}
\newcommand{\hC}{\wh{\mathbb{C}}}
\newcommand{\C}{\mathbb{C}}

\newcommand{\B}{\mathbf{B}}
\newcommand{\T}{\mathbf{T}}

\newcommand{\Om} {\Omega}
\newcommand{\vk} {\varkappa}
\newcommand{\x} {\mathbf x}
\renewcommand{\a} {\alpha}

\begin{document}

\title{Polygons and non-starlikeness of Teichm\"{u}ller spaces}
\author{Samuel L. Krushkal}

\begin{abstract}
The problem of starlikeness of Teichm\"{u}ller spaces in Bers' embedding was raised in 1974 and is solved (negatively) for Teichm\"{u}ller spaces of sufficiently large dimensions.

The original proof given by the author relies on the existence of conformally rigid domains established by Thurston.
Later the author found another proof of non-starlikeness of universal Teichm\"{u}ller space based on geometric features of rectilinear polygons.

This paper provides a complete solution of the problem for Teichm\"{u}ller spaces $\T(g,0)$ of closed Riemann surfaces of genus $g \ge 2$.

\end{abstract}

\date{\today\hskip4mm ({\tt PolygNonStar(2).tex})}

\maketitle

\bigskip


\medskip


\bigskip

\markboth{Samuel L. Krushkal}{Polygons and non-starlikeness of Teichm\"{u}ller spaces}
\pagestyle{headings}

\bigskip\bigskip
\centerline{\bf 1. INTRODUCTORY REMARKS}

\bigskip
There is an interesting still unsolved completely question on shape of holomorphic embeddings of Teichm\"{u}ller  spaces stated in the famous collection of problems on Riemann surfaces, Teichm\"{u}ller spaces and Kleinian groups in the book \cite{BK}:

\bigskip\noindent
{\it For an arbitrary finitely or infinitely generated Fuchsian group $\Gamma$ is the Bers embedding of its Teichm\"{u}ller
space $\mathbf T(\Gamma)$ starlike?}

\bigskip
Recall that in this embedding the space $\mathbf T(\Gamma)$ is represented as a bounded domain formed by the Schwarzian derivatives
$$
S_w = \Bigl(\frac{w^{\prime\prime}}{w^\prime}\Bigr)^\prime - \frac{1}{2} \Bigl(\frac{w^{\prime\prime}}{w^\prime}\Bigr)^2
$$
of holomorphic univalent functions $w(z)$ in the lower half-plane $U^* = \{z: \Im z < 0\}$)
(or in the disk) admitting quasiconformal extensions to the Riemann sphere $\wh{\C} = \C \cup \{\iy\}$ compatible with the group $\Gamma$ acting on $U^*$.

The problem was solved (negatively) for Teichm\"{u}ller spaces of sufficiently large dimensions.
It was shown in \cite{Kr1} that universal Teichm\"{u}ller space $\mathbf T = \mathbf T(\mathbf 1)$ has
points which cannot be joined to a distinguished point even by curves of a considerably general form, in particular, by polygonal lines with the same finite number of rectilinear segments. The proof relies on the existence of conformally rigid domains established by Thurston in \cite{Th} (see also \cite{As}).

This implies, in particular, that universal Teichm\"{u}ller space is not starlike with
respect to any of its points, and there exist points $\varphi \in \mathbf T$ for which the line interval
$\{t \varphi: 0 < t < 1\}$ contains the points from $\mathbf B \setminus \mathbf S$, where
$\mathbf B = \mathbf B(U^*)$ is the Banach space of hyperbolically bounded holomorphic functions in the
half-plane $U^*$ with norm
 \be\label{1}
\|\vp\|_{\mathbf B} = 4 \sup_{U^*} y^2 |\vp(z)|
\end{equation}
and $\mathbf S$ denotes the set of all Schwarzian derivatives of univalent functions on $U^*$. These points correspond to holomorphic functions on $U^*$ which are only locally univalent.

On this way, it was established in \cite{Kr3} that also all finite dimensional Teichm\"{u}ller spaces $\mathbf T(\Gamma)$
of high enough dimensions $n \ge n_0$ are not starlike.

Toki \cite{To} extended the result on the non-starlikeness of the space $\mathbf T$ to Teichm\"{u}ller spaces of Riemann surfaces that contain hyperbolic disks of arbitrary large radius, in particular, for the
spaces corresponding to Fuchsian groups of second kind. The crucial point in the proof of \cite{To} is the same as in \cite{Kr1}.

Later the author found much simpler constructive proof of non-starlikeness of the space $\T$ based on geometric features of rectilinear polygons representing explicitly the functions which violate this property. This is given in \cite{Kr5}.

The non-starlikeness causes obstructions to some problems in the Teichm\"{u}ller space theory and its applications to geometric complex analysis.

\bigskip\bigskip
\centerline{\bf 2. NEW NON-STARLIKENESS RESULTS}

\bigskip
In this paper, we establish by modifying the approach of \cite{Kr5} the complete solution of this problem for  Teichm\"{u}ller spaces $\T(g, 0)$ of closed (hyperbolic) Riemann surfaces of genus $g \ge 2$. We prove the following

\bigskip\noindent
{\bf Theorem 1}. {\it Every space $\T(g, 0)$ with $g \ge 2$ is not starlike with respect to its origin in Bers imbedding}.

\bigskip
As a simple consequence, this theorem implies a similar result for one space of planar surfaces with punctures:

\bigskip\noindent
{\bf Theorem 2}. {\it The space $\T(0, 6)$ of spheres with six punctures is not starlike with respect to the origin in Bers imbedding}.

\bigskip
This theorem follows by using the natural isomorphism between the spaces $\T(0, 6)$ and $\T(2,0)$ regarding the points of the second  space as hyperelliptic surfaces of genus two branched over the sphere $\hC$.

So the question remains open only for the finite dimensional spaces $\T(g, n)$ corresponding to closed Riemann surfaces of genus $g$ with $n \ge 1$ punctures such that
$$
2 < 3g - 3 + n < n_0,
$$
with $n_0$ indicated above.

In the last section we present some generalization of these theorems.

\bigskip\bigskip
\centerline{\bf 3. PROOF OF THEOREM 1}

\bigskip
We accomplish the proof in several steps starting with  some preliminary constructions.
We shall represent the points of $\T(\Gamma)$ by Fuchsian groups $\Gamma$ acting discontinuously on
the disks $\D = \{|z| < 1\}$ and $\D^* = \{z \in \hC: |z| > 1\}$. The hyperbolic metric on $\D$ of
curvature $- 4$ has the differential form
$$
ds = |dz|/(|z|^2 - 1)^2,
$$
and the corresponding space $\B = \B(\D)$ of hyperbolically bounded functions$\vp$ on $\D$  has instead of (1) the norm
$$
\|\vp\|_{\B} = \sup_\D (1 - |z|^2)^2 |\vp(z)|.
$$
We shall also use the corresponding space of holomophic functions in the disk $\D^*$; these functions
satisfy $\vp(z) = O(z^{-4})$ near $z = \iy$.

\bigskip\noindent
$\mathbf{1^0}$. \ First of all, observe that for any Fuchsian group $\Gamma$ of acting on the disk $\D$,
the universal covering $\D \to \D/\Gamma$ naturally generates a canonical embedding of the corresponding Teichm\"{u}ller space $\T(\Gamma)$ into the universal space $\T$. All our groups $\Gamma$ are of the first kind, that means the unit circle $\mathbb S^1 = \partial \D$ is the limit set of $\Gamma$.

Accordingly, we have on each space $\T(\Gamma)$ two Teichm\"{u}ller distances generated by quasiconformal maps: the restricted metric $\wt \tau_{\T(\Gamma)} = \tau_\T|\T(\Gamma)$ and the intrinsic metric $\tau_{\T(\Gamma)}$ generated by maps compatible with the group $\Gamma$.

These metrics are equivalent: clearly, $\wt \tau_{\T(\Gamma)} \le \tau_{\T(\Gamma)}$; some explicit upper bounds for the ratio $\tau_{\T(\Gamma)}/\wt \tau_{\T(\Gamma)}$ are given in \cite{Le}, Ch. V.

Below we shall use on the space $\T(\Gamma)$ the distance $\wt \tau_{\T(\Gamma)}$ coming from the ambient space $\T$  and write, if needed, $(\T(\Gamma), \wt \tau)$.

As is well known, every space $\T(\Gamma)$ of dimension greater than $1$, does not be equivalent holomorphically to a Banach ball. Since
$$
\T(\Gamma) = \T \bigcap \B(\Gamma),
$$
where $\B(\Gamma)$ is the $(3n - 3)$-dimensional linear subspace of $\B$ formed by the $\Gamma$-automorphic forms on $\D$ of weight $- 4$ (quadratic $\Gamma$-differentials), and $\T(\Gamma)$ contains the ball
$\{\|\vp\|_{\B(\Gamma)} < 2\}$, there are the points $\vp \in \T(\Gamma)$ with $\|\vp\|_{\B(\Gamma)} > 2$
(their collection is an open subset); cf. \cite{Ea}, \cite{Le}.

\bigskip\noindent
$\mathbf{2^0}$. The following important lemma is one of the underlying facts in the proof of Theorem 1.

\bigskip\noindent
{\bf Lemma 1}. {\it For any rational function $r_n$ with poles of order two on the boundary circle $\mathbb S^1$ of the form
  \be\label{2}
r_n(z) = \sum\limits_1^n \fc{c_j}{(z - a_j)^2} + \sum\limits_1^n \fc{c_j^\prime}{z - a_j}
\end{equation}
satisfying $\sum\limits_1^n |c_j| > 0$, we have the equality
 \be\label{3}
\|r_n\|_\B = \limsup\limits_{|z| \to 1} (1 - |z|^2)^2 \ |r_n(z)|.
\end{equation}
}

This lemma is a special case of a general deep fact established in \cite{Kr6} for arbitrary quasidisks $D \subset \hC$.

It follows from (3) that there is a boundary point $z_0 \in \mathbb S^1$ at which the maximal value of $(1 - |z|^2)^2 |r_n(z)|$ on the closed disk $\ov \D$ is attained.

Another important fact, which also will be applied in the proof of Theorem 1, is that the Schwarzians of  conformal mapping functions of circular polygons are of the form (2).

\bigskip\noindent
$\mathbf{3^0}$. We proceed  to the proof of Theorem 1 and uniformize the base point $X_0$ of the space $\T(g) = \T(X_0)$ by a Fuchsian group $\Gamma_0$ acting discontinuously  on the disks $\D^*$, in other wrds, consider this point as the quotient space $X_0 = \D^*/\Gamma_0$.

As was mentioned above, there is a surface $X \in \T(X_0)$ corresponding to a point $\vp \in \B(\Gamma_0)$ with
 \be\label{4}
\|\vp\|_{\B(\Gamma_0)} > 2.
\end{equation}
This surface is represented (uniformized) by a quasifuchsian group
$\Gamma_\mu = (f^\mu)^{-1} \Gamma_0 f^\mu$, where $f^\mu(z)$ is a quasiconformal automorphism of $\hC$
conformal on $\D^*$ having the Schwarzian $S_{f^\mu}(z) = \vp(z), \ z \in \D^*$ (the group $\Gamma_\mu$
depends only on $\vp$).

We take Ford's fundamental polygon $P(\Gamma_\mu)$ of this group in domain $D_\mu^* = f^\mu(\D^*)$.
which is the intersection of the exterior domains of all isometric circles of this group:
$$
P(\Gamma_\mu) = \{z \in D_\mu^*: \ \sup \ |\gamma^\prime(z)| < 1\},
$$
taking the supremum over all elements of the group $\Gamma_\mu$ different from the identity.
The sides of this polygon are pairwise $\Gamma_\mu$-equivalent. Its inverse image $(f^\mu)^{-1} P(\Gamma_\mu)$ is
the fundamental domain of the initial group $\Gamma_0$ in $\D^*$.
Since $X$ is a closed surface, the interior and exterior angles of both polygons $P(\Gamma_\mu)$
and $(f^\mu)^{-1} P(\Gamma_\mu)$ are different from $0, \ \pi$ and $2 \pi$.

Now, denoting the elements of $\Gamma_\mu$ by $\gamma_1 = \mathbf 1, \ \gamma_2, \gamma_3, \ \dots $,
numerated so that two successive $\gamma_j, \gamma_{j+1}$ determine the adjacent polygons
$ \gamma_j P(\Gamma_\mu), \  \gamma_{j+1} P(\Gamma_\mu)$ with a common side
consider the increasing connected unions
$$
\Om_1 = P(\Gamma_\mu), \quad \Om_2 = \gamma_1 P(\Gamma_\mu) \bigcup \gamma_2 P(\Gamma_\mu),
\quad \Om_3 = \gamma_1 P(\Gamma_\mu) \bigcup \gamma_2 P(\Gamma_\mu) \bigcup\gamma_3 P(\Gamma_\mu),        \dots
$$
exhausting the domain $D_\mu^*$. All domains $\Om_j$ represent the circular polygons with nonzero angles containing inside the infinite point.

Consider the corresponding conformal maps $F_j$ of the disk $\D^*$ onto the polygons $\Om_j$
and their Schwarzians $S_{F_j}(z)$. Due to the Ahlfors-Weill theorem \cite{AW}, all functions $t S_{F_j(z)}$ with $\|t S_{F_j(z)}\|_\B < 2$ are the Schwarzian derivatives of univalent function $W_t(z)$
on $\D^*$ having quasiconformal extensions to $\D$ with harmonic Beltrami coefficients
 \be\label{5}
\nu_{F_j}(z) = - \fc{1}{2} (1 - |z|^2)^2 \ S_{F_j}(1/\ov z), \quad |z| <  1.
\end{equation}

Lemma 1 (actually, the equality (3)) allows one to apply to any $F_j$ rather standard arguments already applied in \cite{Kr5}, \cite{Kr7}, which provide that the harmonic Beltrami coefficients (4) are extremal in their equivalence classes $[F^{\nu_{S_j}}]$ (i.e., among quasiconformal extensions of $F_j$ onto $\D$). Note that {\it the extremality here is doing with respect to the metric $\wt \tau_{\T(G)}$ }
and means that the extremal rays are geodesic in this metric.

Fix a subarc $\gamma \subset \mathbb S^1$ and take the conformal map $z = \chi(\zeta)$ of the half-strip
$$
\Pi_{+} = \{\zeta = \xi + i \eta : \ \xi > 0, \ 0 < \eta < 1\}
$$
onto $\D^*$ such that $\chi_1^{-1}(z_0) = \iy$ and the pre-images of the endpoints of the arc $\gamma$ are the points $0$ and $1$.
Then the coefficient $\nu_{t S_{F_j}}$ is pulled back to Beltrami coefficient
$$
t \mu_j(\z) := \chi_1^* \nu_{t S_{F_j}}(z) = (\nu_{t S_{F_j}} \circ \chi)(\z) \ \ov{\chi^\prime(\z)}/\chi^\prime(\z)
$$
on $\Pi_{+}$. We show that either from these coefficients has an essential boundary point
at infinity, and  a degenerating sequence quadratic differentials at this point;
this sequence consists of the functions
$$
\omega_m(\z) = \frac{1}{m} e^{- \z/m}, \quad \z \in \Pi_{+}, \ \ m = 1, 2, ... \ .
$$
All these functions belong to $A_1^2(\Pi_{+})$ and satisfy $\omega_m(\zeta) \to 0$ uniformly on $\Pi_{+} \bigcap \{|\z| < M\}$ for any $M < \iy$; in addition,
$\big\vert \iint_{\Pi_{+}} \omega_m(\z) d\xi d\eta \big\vert = 1 - O(1/m)$.

Our goal is to show that for any $j$,
  \be\label{6}
\lim\limits_{m\to \iy} \left | \iint\limits_{\Pi_{+}} \mu_j(\zeta) \omega_m(\zeta) d \xi d \eta \right | = \|\mu_j\|_\iy.
\end{equation}
This equality follows from the relations
 \be\label{7}
\int\limits_0^\iy \fc{\partial \mu_j(\xi + i \eta)}{\partial \xi} e^{- \xi/m} d \xi
= \fc{1}{m} \int\limits_0^\iy \mu_j(\xi + i \eta) e^{- \xi/m} d \xi - \mu_j(i \eta)
= \fc{1}{m} \int\limits_0^\iy \mu_j(\xi + i \eta) e^{- \xi/m} d \xi.
\end{equation}
The last integral in (7) represents the values of the Laplace transform
$$
\mathcal L \nu = \int\limits_0^\iy \nu(t) e^{-st} dt
$$
of $\nu(\xi + i \eta)$ in the points $s = 1/m$. Taking into account the Lebesgue theorem on dominated convergence  and the properties of the extremal Beltrami coefficient $\mu_j$ (in the metric  $\wt \tau_{\T(\Gamma)}$), one derives
$$
\lim\limits_{m \to \iy} \fc{1}{m} \Bigl | \int\limits_0^\iy
\mu_j(\xi + i \eta) e^{- \xi/m} d \xi \Bigr | = |\mu_j(\iy) - \nu(i \eta)| = |\mu_j(\iy)|
$$
(where $\mu_j(\iy) = \lim\limits_{\xi \to - \iy} \mu_j(\xi + i \eta)$). Since
$$
\lim\limits_{m\to \iy} \left | \iint\limits_{\Pi_{+}} \nu_{t S_{F_j}}(\zeta) \omega_m(\zeta) d \xi d \eta \right |
= \lim\limits_{m\to \iy} \left | \iint\limits_D \mu_j(z) \psi_m(z) dx dy \right | = \|\nu_{t S_{F_j}}\|_\iy,
$$
with
 \be\label{8}
\psi_m = (\omega_m \circ \chi^{-1})(\chi^\prime)^{-2}, \quad m = 1,2, \dots,
\end{equation}
we have
 \be\label{9}
\lim\limits_{m\to \iy}
\Bigl | \iint\limits_D \mu_j(z) \psi_m(z) dx dy \Bigr |
= \lim\limits_{m\to \iy} \Bigl | \iint\limits_{\Pi_{+}} \nu_{t S_{F_j}}(\zeta) \omega_m(\zeta) d \xi d \eta \Bigr | = |\mu_j(\iy)|,
\end{equation}
which implies (6). One also obtains that the sequence (8) is degenerated for $\nu_{t S_{F_j}}$.
Consequently, $z_0$ is a substantial point of $\nu_{t S_{F_j}}$ and this Beltrami coefficient is extremal in its equivalence class (hence, this class cannot contain a Teichm\"{u}ller Beltrami coefficient).

The equality (6) implies the extremality of all Beltrami coefficients (4) for the metric $\wt \tau_{\T(G)}$  and also  yields that all maps $F_j$ have equal Teichm\"{u}ller and Grunsky norms.

\bigskip\noindent
$\bf{4^0}$. \ By the general Carath\'{e}odory theorem on convergence of conformal maps, the limit function
 \be\label{10}
F(z) = \lim\limits_{j \to \iy} F_j(z)
\end{equation}
maps conformally the disk $\D^*$ onto the exterior component $D_\mu^*$ of domain of discontinuity of the group $\Gamma_\mu$.
Note that the convergence in (10) is uniform on the disk $\D^*$ (in the spherical metric on $\hC$).

Our goal now is to establish that one of the extremal qusiconformal extensions $\wh F^\mu$ of this map inherits the properties of the functions $F_j$. This is the crucial step in the proof of Theorem 1.

It follows from (10) that the Beltrami coefficients (5) are convergent on the disk $\D$ to the harmonic
Beltrami coefficient $\nu_{S_F}$ of function $F$; besides,
 \be\label{11}
\|S_F\|_\B \le \lim\limits_{j \to \iy} \|S_{F_j}\|_\B
\end{equation}
(this inequality also follows from the extremality of $\nu_{tS_{F_j}}$ and general properties of quasiconformal maps).

We show that in fact we have in (11) the equality, which yields that for admissible $|t|$, the harmonic
Beltrami coefficient $\nu_{t S_F}$ is extremal (in the distance $\wt \tau_\T(\Gamma)$) in its equivalence class.

First observe that normalising the quasiconformal extensions $\wh F(z)$ of univalent functions $F(z) = z + b_0 + b_1 z^{-1}+ \dots$ on $\D^*$ to $\hC$ additionally by $\wh F(0) = 0$, one obtains that their
Teichm\"{u}ller norm $k(F)$ is plurisubharmonic on the universal Teichm\"{u}ller space $\T$ (hence,
also  on $\T(\Gamma_\mu)$ with distance $\wt \tau_\T(\Gamma)$).

Now we pass to the homotopy functions
$$
F_s(z) = s F(z/s) = z + b_0 s + b_1 s^2 z^{-1}+\dots
$$
with $|s| \le 1$. Then $k(F_s)$ is circularly symmetric with respect to $c$, hence continuous on $[0, 1]$.

But for any $s \in (0, 1)$, we have
$$
\|S_{F_{j,s}} - S_{F_s}\|_\B \to 0 \quad \text{as} \ \ j \to \iy,
$$
and the properties of extremal quasiconformal maps imply
 \be\label{12}
\lim\limits_{s \to 1} k(F_{j,s}) = k(F_j) = \|\nu_{S_F}\|_\iy
\end{equation}
(the last equality in (12) is valid only for $\|S_F\|_\B < 2$).
Together with (4) and (7), this implies the equality in (9).

It follows from (12) that all harmonic Beltrami coefficients $\nu_{t S_F}$ with norm less than $1$ {\bf must be extremal in their classes}, and therefore, as $\|\nu_{t S_F}\|_\B \to 1$, the corresponding value  $t_0$ must define a boundary point of both spaces $\T(\Gamma_\mu)$ and $\T$ (i.e., $t S_F \to t_0 S_F \in \partial \ T(\Gamma_\mu)$ as $t \to t_0$).

Taking also into account the Abikoff-Bers-Zhuravlev theorem, which states that any the domain $\ T(\Gamma)$ representing the space $\T(\Gamma)$ has a common boundary with its complementary domain in the space $\B(\Gamma)$ (see \cite{Ab}, \cite{Be}, \cite{Zh}), one derives that there are points on the ray ${r S_F: \ r > 0}$ in $\T(\Gamma_\mu)$, which must lie in the exterior of $\T(\Gamma_\mu)$.
These points are placed on the interval between $r_0 S_F $ and $r_{*} S_F$, where the second point correspond to the surface $D_\mu^*/\Gamma_\mu$ obeying the assumption (4).
This completes the proof of the theorem.

\bigskip\bigskip
\centerline{\bf 4. ADDITIONAL REMARKS}

\bigskip\noindent
{\bf 1. Alternate proof of relarions (12)}.
In view of importance, we also present another proof of the relations (12). It involves the Grunsky coefficients of univalent function.

Recall that these coefficients are defined from the expansion
$$
\log \fc{f(z) - f(\z)}{z - \z} = - \sum\limits_{m,n=1}^\iy \a_{mn} z^{-m} \z^{-n}, \quad (z, \z) \in (\D^*)^2,
$$
and by the Grunsky theorem \cite{Gr}, a holomorphic function $f(z) = z + \const + O(z^{-1})$
in a neighborhood of $z = \iy$ can be extended to a univalent holomorphic function on $\D^*$ if and only
the coefficients $\a_{m n}$ satisfy the inequality
$$
\Big\vert \sum\limits_{m,n=1}^\iy \ \sqrt{m n} \ \a_{m n} x_m x_n \Big\vert \le 1,
$$
for any sequence $\mathbf x = (x_n)$ from the unit sphere $S(l^2)$ of the Hilbert space $l^2$ with norm $\|\x\| = (\sum\limits_1^\iy |x_n|^2)^{1/2}$ (here the principal branch of the logarithmic function is chosen). The quantity
$$
\vk(f) = \sup \Big\{\Big\vert \sum\limits_{m,n=1}^\iy \ \sqrt{mn} \ \a_{m n} x_m x_n \Big\vert: \
\mathbf x = (x_n) \in S(l^2) \Big\} \le 1
$$
is called the Grunsky norm of $f$. Also this norm is a plurisubharmonic function on the space $\T$.
It is majorated by the Teichm\"{u}ller norm $k(F^\mu)$, and the equality $\vk(F^\mu)= k(F^\mu)$
is valid if and only if the extremal Beltrami coefficient $\mu_0$  in the equivalence class of $F^\mu$ satisfies
$$
\|\mu_0\|_\iy = \sup_{\psi \in A_1^2(\D),\|\psi\|_{A_1} =1}
\Big\vert \iint\limits_\D \mu_0(z) \psi(z) dx dy \Big\vert \quad (z = x + iy \in \D),
$$
where $A_1(\D^*)$ is the space of integrable holomorphic quadratic differentials $\psi(z) dz^2$ on
$\D^*$ and
$$
A_1^2(\D) = \{\psi = \om^2 \in A_1(\D): \ \om \ \ \text{holomorphic \ in} \ \ \D\}.
$$

The equality $\vk(F_j) = k(F_j)$ following from (7) implies the similar relation for the homotopy
maps $F_{j,s}$.
Another important deep fact (first established by K\"{u}hnau \cite{Ku}, but also following from harmonicity
of $\vk(F_s^\mu)$ in $t \in \D$) is that along the homotopy disk,
$$
\lim\limits_{s \to 1} \vk(F_s^{\mu}) = \vk(F^\mu).
$$
Combining all this with the basic inequality  $\vk(F) \le k(F)$, one obtains the equalities (12).

\bigskip\noindent
{\bf 2}. As a consequence from the second proof, one obtains the following assertion.

\bigskip\noindent
{\bf Proposition 1}. {\it For all Beltrami coefficients $\mu$ compatible with the groups
$\Gamma$  uniformizing the closed Riemann surfaces, we have the equality $\vk(F^\mu) = k(F^\mu)$
(in the Teichm\"{u}ller distance on $\T$).  }

\bigskip
This result has its intrinsic interest, because the points $S_F \in \T$ with $\vk(F) < k(F)$ are dense in
the space $\T$.

\bigskip\noindent
{\bf 3. Generalization of Theorems 1 and 2}. The arguments implied in the proof of Theorem 1 provide,
in fact, more general  result.

\bigskip\noindent
{\bf Theorem 3}. {\it Let $\Gamma$ be a Fuchsian group such that $X_0 = \D/\Gamma$ is a closed Riemann surface of genus $g \ge 2$, and let $G$ be a subdomain of $\T(\Gamma) \subset \B(\Gamma)$ satisfying: 

(a) $G$ contains the origin $\mathbf 0$ of $\B(\Gamma)$ and a point $\vp_0$ with $\|\vp_0\|_\B > 2$;

(b) the domain $G$ and its complementary domain $G^*$ in $\B(\Gamma)$ have a common boundary
$\partial G$.

Then $G$ is not starlike with respect to the origin.   }

Theorem of such type is also valid for appropriate subdomains of the space $\T(0, 6)$.

\medskip
{\small\em{ \leftline{Department of Mathematics, Bar-Ilan University, Ramat-Gan, Israel}
\leftline{and Department of Mathematics, University of Virginia, Charlottesville, VA 22904-4137, USA}}}

\end{document}